\documentclass{iise}

\usepackage{multirow}
\usepackage{optidef}
\RequirePackage{subcaption}
\usepackage{tabularx}

\conference{Proceedings of the IISE Annual Conference \& Expo 2024 \\
A. Brown Greer, C. Contardo, J.-M. Frayret, eds.}

\title{\titlesize Pediatric vaccine tender scheduling in low- and middle-income countries}

\author{
Nicholas Uhorchak$^{1}$,
Ruben A. Proano$^{2}$,
Sandra Eksioglu$^{1}$,
Fatih Cengil$^{1}$,
Burak Eksioglu$^{1}$
\\
$^{1}$ Department of Industrial Engineering,\\University of Arkansas, Fayetteville, AR 72701\\
\vspace{0.3cm}
$^{2}$ Department of Industrial and Systems Engineering,\\Rochester Institute of Technology, Rochester, NY 14623\\
}

\authorlist{Uhorchak, Proano, Eksioglu, Cengil, and Eksioglu}
\abstractID{7176}

\begin{document}
\maketitle

\begin{abstract}
{\small Effective and efficient scheduling of vaccine distribution can significantly impact vaccine uptake, which is critical to controlling the spread of infectious diseases. Ineffective scheduling can lead to waste, delays, and low vaccine coverage, potentially weakening the efforts to protect the public. Organizations such as UNICEF (United Nations Children’s Fund), PAHO (Pan American Health Organization), and GAVI (Gavi, the Vaccine Alliance) coordinate vaccine tenders to ensure that enough supply is available on the international market at the lowest possible prices. Scheduling vaccine tenders over a planning horizon in a way that is equitable, efficient, and accessible is a complex problem that involves trade-offs between multiple objectives while ensuring that vaccine availability, demand, and logistical constraints are met. The current method for scheduling tenders is generally reactive and over short planning horizons. Vaccine tenders are scheduled when supply is insufficient to cover demand. We propose an optimization model to dynamically and proactively generate vaccine tender schedules over long planning horizons. This model helps us address the following research questions: What should the optimal sequencing and scheduling of vaccine tenders be to enhance affordability and profit over long time horizons? What is the optimal tender procurement schedule for single/multiple antigen scenarios? We use several real-life data sources to validate the model and address our research questions. Results from our analysis show when to schedule vaccine tenders, what volumes manufacturers should commit to, and the optimal tender lengths to satisfy demand. We show that vaccine tenders tend towards maximum lengths, generally converge over long time horizons, and are robust to changes in varying conditions.}
\end{abstract}

\section*{Keywords}
Integer programming, vaccine tender, hypothetically coordinated vaccine market

\section{Introduction}
Pediatric vaccines prevent an estimated "3.5-5 million deaths from diseases like diphtheria and tetanus" annually\cite{WHO2024}. Immunization campaigns prevent diseases through the administration of vaccines that induce the production of antibodies against specific disease-causing antigens. 

To meet market demand, vaccine tenders for different antigens are scheduled to ensure adequate supply levels. Low- and low-middle-income countries (LIC, and LMIC, respectively) primarily procure vaccines in the UNICEF market, and depending on their gross national income per capita and the strength of their healthcare systems, can receive financial support from GAVI to procure vaccines. UNICEF and PAHO act as coordinating entities, working with countries to aggregate vaccine demand information. UNICEF coordinates vaccine tenders with vaccine producers to meet the demand of LIC and LMICs. {UNICEF plans tenders 2 to 3 years in advance. However, its scheduling is reactive, primarily considering short-term demand needs, having limited visibility of other factors, such as country-specific immunization plans, and manufacturing capacity producers share among multiple products in their portfolio.}

In a tender, manufacturers commit to supply specific vaccine volumes to the UNICEF market at an agreed price for a given number of years (i.e., tender length). We propose an analytical approach to determine when to schedule vaccine tenders, what volumes manufacturers should commit to, and the optimal tender lengths to satisfy the long-term needs of LIC and LMIC in a hypothetically coordinated vaccine market (HCVM). An HCVM involves vaccine producers, LIC and LMIC grouped in market segments, and altruistic coordinating entities like UNICEF or PAHO. These entities aim to meet the vaccine demands at affordable prices while ensuring producers' returns. Entities like UNICEF and PAHO manage procurement, logistics, and financial support, ensuring a steady vaccine supply that meets the countries' immunization needs within their economic constraints. To this end, we develop a mixed-integer program that identifies vaccine tenders' optimal timing and duration over long planning horizons. Our findings show that {for varying conditions, it is possible to schedule vaccine tenders longer than the current practice of 2-3 years and provide insight into possible vaccine manufacturer commitment.}

Various mathematical models exist for scheduling pediatric vaccines \cite{YuVax}. However, these efforts do not focus on the long-term scheduling of vaccine procurement tenders. Additional studies thoroughly review vaccine procurement and highlight the associated processes and policy implications \cite{Gianfredi2021}. Long-term procurement scheduling is most commonly found in energy procurement models \cite{Zhang2018}.

Proano et al. \cite{Proano2012} introduce a deterministic model for the HCVM to explore if it is possible to find vaccine prices that can be affordable and profitable for a coordinated one-year multi-vaccine tender. In a subsequent study, Proano and Mosquera \cite{Proano2021} establish that implementing tier-pricing with more segments within the GAVI and PAHO markets could extend the opportunities to have affordable vaccines. Alves-Maciel and Proano \cite{Maciel2021, Maciel2022-1} further extend HCVM to study the effect of market coordination under multiple non-cooperative coordinating entities and when tenders rely on lump-sum payments for vaccine product baskets.

The existing research on the HCVM has focused on the impact of market coordination on vaccine affordability and profit for a single multi-product short-term tender. In contrast, this study focuses on long-term tender scheduling for several antigens. The objective is to enhance availability and profitability by enabling coordinating entities to implement long-term tender schedules. The model considers that the available capacity of manufacturers, and potential capacity expansions, are shared among different antigens. 

\section{Model Description}

\newcolumntype{A}{p{0.59\textwidth}}
\newcolumntype{B}{p{0.39\textwidth}}

\textbf{Sets}: \vspace{-3mm}
\begin{table}[h]
\small
\label{sets}
    \begin{tabular}{AB}
    $\mathcal{V}$: set of vaccines (from the WHO pre-qualified vaccine list \cite{WHO_PQ}), &  $\mathcal{A}$: set of antigens, \\
    $\mathcal{V}_a$: set of vaccines containing antigen $a$, & $\mathcal{P}$: set of all producers, \\
    $\mathcal{A}_p$: subset of antigens contained in producer $p$'s portfolio, & $\mathcal{P}_v$: set of producers of vaccine $v$, \\
    $\mathcal{T}$: set of time periods in the planning horizon, & $\mathcal{A}_v$: subset of antigens in vaccine $v$, \\
    \multicolumn{2}{l}{$\mathcal{T}_t'$: set of periods when a tender ends \big(since we allow tenders of at most 5 years, $\mathcal{T}_t' = \{t, t+1, ..., min\{t+4,|\mathcal{T}|\}\}$\big).} 
    \end{tabular}
\end{table}

\textbf{Parameters}: \vspace{-3mm}
\begin{table}[h]
\small
\label{parameters}
    \begin{tabular}{AB}
    $l_{vp}$: the annualized return on investment required by producer $p$ for vaccine $v$, & ${d}_{at}$: demand for antigen $a$ for period $t$, \\
    $f_{vpt}$: the production set-up cost of producer $p$ for vaccine $v$ in period $t$, & ${s}_{pt}$:  production capacity of producer $p$ for period $t$, \\
    $h_{v}$: the annual holding cost for vaccine $v$ as a proportion of its price, & $g_t$: the set-up cost for initiating a tender in period $t$, \\
    $r_{vpt}$: reservation price of vaccine $v$ produced by $p$ in period $t$, & $\Gamma_p$: capacity extension cost for producer $p$, \\
    $\bar{r}_{vt}$: average price of vaccine $v$ across all producers in period $t$, & $\kappa$: allowable capacity extension rate, \\
    $S_{a0}$:initial number of children unvaccinated with antigen $a$, & $I_{v0}$: initial inventory for vaccine $v$, \\
    $\beta$ : penalty per unit of shortage of the committed amount, & $\delta$: the discount factor\\
    \hspace{0.17in}i.e., the unvaccinated children penalty. 
    \end{tabular}
\end{table}

\newcolumntype{C}{p{0.95\textwidth}}

\textbf{Decision Variables}: \vspace{-3mm}
\begin{table}[h]
\small
\label{decisions}
    \begin{tabular}{C}
    $F_{at\tau}$ is a binary variable indicating whether a tender covers the demand for antigen $a$ for periods $t$ through $\tau$, \\
    $Q_{vpt\tau}$ is the procurement commitment for vaccine $v$ by producer $p$ for periods $t$ through $\tau$, \\
    $X_{vpt}$ represents doses of vaccine $v$ delivered by producer $p$ in period $t$, \\
    $Y_{pt}$ is a binary variable which specifies whether or not manufacturer $p$ produces in period $t$, \\
    $W_{pt\tau}$ indicates if producer $p$ has made commitments for periods $t$ through $\tau$, \\
    $L_{pt}$ is a binary variable indicating if a capacity extension decision has been made, \\
    $I_{vt}$ is the stock level for vaccine $v$ at the beginning of period $t$, \\
    $V_{vt}$ is the number of children vaccinated with vaccine $v$ in period $t$, \\
    $S_{at}$ is the number of children unvaccinated for antigen $a$ in period $t$.
    \end{tabular}
\end{table}

\textbf{Formulation}: \\
Using the notation provided above, the problem is formulated as follows:

{\small
\begin{mini!}
{}{\sum_{t\in \mathcal{T}} \delta^t \left(\sum_{a \in \mathcal{A}}  \sum_{\tau \in \mathcal{T}_t'} g_t F_{at\tau} + \sum_{v \in \mathcal{V}} \sum_{p \in \mathcal{P}_{v}} {r_{vpt}} X_{vpt} + \sum_{a \in \mathcal{A}}   \beta S_{at} + \sum_{v \in \mathcal{V}} h_v\bar{r}_{vt}I_{vt} + \sum_{p \in \mathcal{P}} \Gamma_p L_{pt}\right)}{}{} \label{obj}
\addConstraint{(\tau-t+1) F_{at\tau}}{\leq \sum_{p\in \mathcal{P}_a}\sum_{l=t}^{\tau}Y_{pl},}{\forall a\in \mathcal{A}; t \in \mathcal{T}; \tau \in \mathcal{T}_t'} \label{eq:tender}
\addConstraint{F_{at\tau} + F_{at\tau'}}{\leq  1,}{\forall a\in \mathcal{A}; t \in \mathcal{T}; \tau, \tau' \in \mathcal{T'}; \tau \neq \tau'} \label{eq:overlap}
\addConstraint{\sum_{t=max\{1,l-4\}}^{l} \sum_{\tau \in \mathcal{T}_t'} F_{at\tau}}{\geq  1,}{\forall a\in \mathcal{A}; l \in T} \label{eq:1ab-newa}
\addConstraint{\sum_{a\in \mathcal{A}_p} F_{at\tau}}{\geq W_{pt\tau},}{\forall  p \in \mathcal{P}; t \in \mathcal{T}; \tau \in \mathcal{T}_t'} \label{eq:2ab-newa} 
\addConstraint{\sum_{a\in \mathcal{A}_p} F_{at\tau}}{\leq |\mathcal{A}_p| W_{pt\tau},}{\forall  p \in \mathcal{P}; t \in \mathcal{T}; \tau \in \mathcal{T}_t'} \label{eq:2ab-newb}
\addConstraint{\sum_{v \in \mathcal{V}_p} Q_{vpt\tau}}{\leq W_{pt\tau} \left(\sum_{l=t}^{\tau}(s_{pl}+\sum_{k=1}^{l}s_{pk}\kappa L_{pk})\right),}{\forall  p \in \mathcal{P}; t \in \mathcal{T}; \tau \in \mathcal{T}_t'} \label{eq:3_a-newa}
\addConstraint{Q_{vpt\tau}}{\geq  W_{pt\tau} \sum_{l=t}^{\tau} X_{vpl},}{\forall v\in \mathcal{V}; p \in \mathcal{P}_v; t \in \mathcal{T}; \tau \in \mathcal{T}_t'} \label{eq:4_0_newa}
\addConstraint{\sum_{v\in \mathcal{V}_p}X_{vpt}}{\leq  Y_{pt} \left(s_{pt}+\sum_{l=1}^{t}s_{pl}\kappa L_{pl}\right),}{\forall p \in \mathcal{P}; t \in \mathcal{T}} \label{eq:4_0-newa2}
\addConstraint{I_{vt-1} +  \sum_{p\in\mathcal{P}_v} X_{vpt}}{= V_{vt} + I_{vt},}{\forall v \in \mathcal{V}; t \in \mathcal{T}} \label{eq:6-newa}
\addConstraint{{d}_{at} - \sum_{v \in \mathcal{V}_a}V_{vt}+ S_{at-1}}{\leq S_{at},}{\forall a\in \mathcal{A};  t \in \mathcal{T}} \label{eq:7-newa}
\addConstraint{\sum_{v\in\mathcal{V}_p} {r_{vpt}} X_{vpt}}{\geq \sum_{v\in\mathcal{V}_p} (1 + l_{vp})f_{vpt}Y_{pt},}{\forall p\in \mathcal{P}; t \in \mathcal{T}} \label{eq:9-newa}
\addConstraint{Q_{vpt\tau}, X_{vpt}, I_{vt}, S_{at}, V_{vt}}{\geq 0,}{\forall a\in \mathcal{A}; v\in \mathcal{V}; p \in \mathcal{P}_v; t \in \mathcal{T}; \tau \in \mathcal{T}_t'} \label{eq:nonnegativity}
\addConstraint{F_{at\tau}, Y_{pt}, W_{pt\tau}, L_{pt}}{\in \{0,1\},}{\forall a\in \mathcal{A}; p \in \mathcal{P}_v; t \in \mathcal{T}; \tau \in \mathcal{T}_t'.} \label{eq:binary}
\end{mini!}
}

The objective function (\ref{obj}) minimizes the total discounted cost, including the setup cost of organizing tenders, vaccine purchasing costs, inventory holding costs, capacity extension costs, and costs of having unvaccinated children. Note that while no true dollar amount can be placed on the long-term costs of unvaccinated children, the penalty is a numerical attempt to prevent unvaccinated children. Constraint (\ref{eq:tender}) ensures that if a tender takes place between years $t$ through $\tau$, at least one producer should supply vaccines during this period. Constraint (\ref{eq:overlap}) ensures that two tenders for the same antigen cannot start at the same period. Constraint (\ref{eq:1ab-newa}) guarantees that at least one tender covers the demand for an antigen in each period. Constraints (\ref{eq:2ab-newa}) and (\ref{eq:2ab-newb}) ensure that a tender is granted to a manufacturer between periods $t$ and $\tau$ if the manufacturer has made commitments for a vaccine for the corresponding antigen ($A_p$). Constraint (\ref{eq:3_a-newa}) assures that a manufacturer's total commitment to the vaccines cannot surpass its production capacity. Note that production capacity in a period might increase based on the capacity extension decision of the manufacturer. Constraint (\ref{eq:4_0_newa}) ensures that the total number of vaccines delivered over given tender periods cannot exceed the vaccine procurement commitment. {Constraint (\ref{eq:4_0-newa2}) ensures that for any given period $t$, the total number of vaccines delivered by a producer does not exceed that producer's production capacity.} Constraint (\ref{eq:6-newa}) ensures inventory balance between consecutive years. Constraint (\ref{eq:7-newa}) keeps track of the number of unvaccinated children. Constraint (\ref{eq:9-newa}) guarantees that each producer makes a profit that fulfills their desired ROI. Constraints (\ref{eq:nonnegativity}) and (\ref{eq:binary}) define the non-negativity and binary variables.

Constraint (\ref{eq:4_0_newa}) includes a non-convex bi-linear term, where a binary variable ($W_{pt\tau}$) and summation of continuous variables ($\sum_{l=t}^{\tau} X_{vpl}$) are multiplied. We linearize (\ref{eq:4_0_newa}) by using McCormick linearization. To create the McCormick envelope, we replace the bi-linear term with a new variable ($\tilde{K}_{vpt\tau}$) and add bounds for this new variable. In the reformulation, constraint (\ref{eq:4_0_newa}) is replaced with the following constraint sets, where the lower and upper bounds for the variables $W_{pt\tau}$ and $\tilde{X}_{vpt\tau}$ are, respectively, defined as $w_{pt\tau}^l$, $w_{pt\tau}^u$, $\tilde{x}_{vpt\tau}^l$, $\tilde{x}_{vpt\tau}^u$. This is a tight relaxation since $W_{pt\tau}$ is a binary variable.

\begin{equation}
\tilde{X}_{vpt\tau} =  \sum_{l=t}^{\tau} X_{vpl}, \hspace{0.3in} \forall v\in \mathcal{V}, p \in \mathcal{P}_v, t \in \mathcal{T}, \tau \in \mathcal{T}_t' \label{McCormick_1}
\end{equation}
\vspace{-1.5mm}
\begin{equation}
Q_{vpt\tau} \geq  K_{vpt\tau}, \hspace{0.3in} \forall v\in \mathcal{V}, p \in \mathcal{P}_v, t \in \mathcal{T}, \tau \in \mathcal{T}_t' \label{McCormick_2}
\end{equation}
\vspace{-1.5mm}
\begin{equation}
K_{vpt\tau} \geq \tilde{x}_{vpt\tau}^u W_{pt\tau} + \tilde{X}_{vpt\tau} w_{pt\tau}^u - \tilde{x}_{vpt\tau}^u w_{pt\tau}^u , \hspace{0.3in} \forall v\in \mathcal{V}, p \in \mathcal{P}_v, t \in \mathcal{T}, \tau \in \mathcal{T}_t' 
\label{McCormick_3}
\end{equation}
\vspace{-1.5mm}
\begin{equation}
K_{vpt\tau} \geq \tilde{x}_{vpt\tau}^l W_{pt\tau} + \tilde{X}_{vpt\tau} w_{pt\tau}^l - \tilde{x}_{vpt\tau}^l w_{pt\tau}^l , \hspace{0.3in} \forall v\in \mathcal{V}, p \in \mathcal{P}_v, t \in \mathcal{T}, \tau \in \mathcal{T}_t' \label{McCormick_4}
\end{equation}
\vspace{-1.5mm}
\begin{equation}
K_{vpt\tau} \leq \tilde{x}_{vpt\tau}^u W_{pt\tau} + \tilde{X}_{vpt\tau} w_{pt\tau}^l - \tilde{x}_{vpt\tau}^u w_{pt\tau}^l , \hspace{0.3in} \forall v\in \mathcal{V}, p \in \mathcal{P}_v, t \in \mathcal{T}, \tau \in \mathcal{T}_t' \label{McCormick_5}
\end{equation}
\vspace{-1.5mm}
\begin{equation}
K_{vpt\tau} \leq \tilde{x}_{vpt\tau}^l W_{pt\tau} + \tilde{X}_{vpt\tau} w_{pt\tau}^u - \tilde{x}_{vpt\tau}^l w_{pt\tau}^u , \hspace{0.3in} \forall v\in \mathcal{V}, p \in \mathcal{P}_v, t \in \mathcal{T}, \tau \in \mathcal{T}_t' \label{McCormick_6}
\end{equation}

\section{Data Collection and Analysis}
The model incorporates synthetically generated data and various datasets from real cases. We focus our analysis on the antigens in Table \ref{tab:vaccines_antigens}. We use the Linksbridge {Vaccine Almanac} datasets, which provide information on vaccine demand in LIC and LMIC in the next 10 years \cite{Linksbridge}. We use the simple exponential smoothing approach to forecast the demand for vaccination beyond 10 years. Production capacities, which often are confidential information, are synthetically estimated using a combination of UNICEF vaccine coverage {estimates}, the average portfolio of producers, and {UNICEF} vaccine demand \cite{UNICEF_supply}. We first aggregate vaccine demand for LIC and LMIC. We then scale the demand based on estimated vaccine coverage. We assume that the production capacity of each manufacturer is proportional to the market share of WHO-prequalified vaccines.
\begin{table}[h]
\centering
\begin{minipage}{\textwidth} 
\begin{tabularx}{\textwidth}{|l|X|} 
\hline
\textbf{Vaccine Name} & \textbf{Antigens} \\
\hline
MMR & Measles, Mumps, Rubella \\
Penta, Hexa & Diphtheria, Tetanus, Pertussis, Haemophilus influenzae type b, Polio, Hepatitis B \\
Rotavirus & Rotavirus \\
HPV & Human papillomavirus \\
PCV & Pneumococcal conjugate \\
\hline
\end{tabularx}
\setlength{\abovecaptionskip}{2pt}
\caption{Vaccines and corresponding antigens}
\label{tab:vaccines_antigens}
\end{minipage}
\end{table}

We used historical data related to the cost of vaccines, which UNICEF publishes \cite{UNICEF_supply}. Since the cost of a vaccine dose differs by manufacturer, we use the corresponding average cost per dose in the model. 
 The reservation price of a vaccine represents the price at which UNICEF purchases vaccines from manufacturers in a specific market. We generate tender schedules for planning horizons ($\mathcal{T}$) ranging from 10 to 30 years. The analysis considers that monovalent and combination vaccines are available to meet vaccination needs. We evaluate the impact of problem parameters on the tender schedule via our numerical analysis.

\section{Discussion of Results}

\begin{figure}[h]
  \centering
  \begin{subfigure}{0.49\textwidth}
    \centering
    \includegraphics[width=\linewidth]{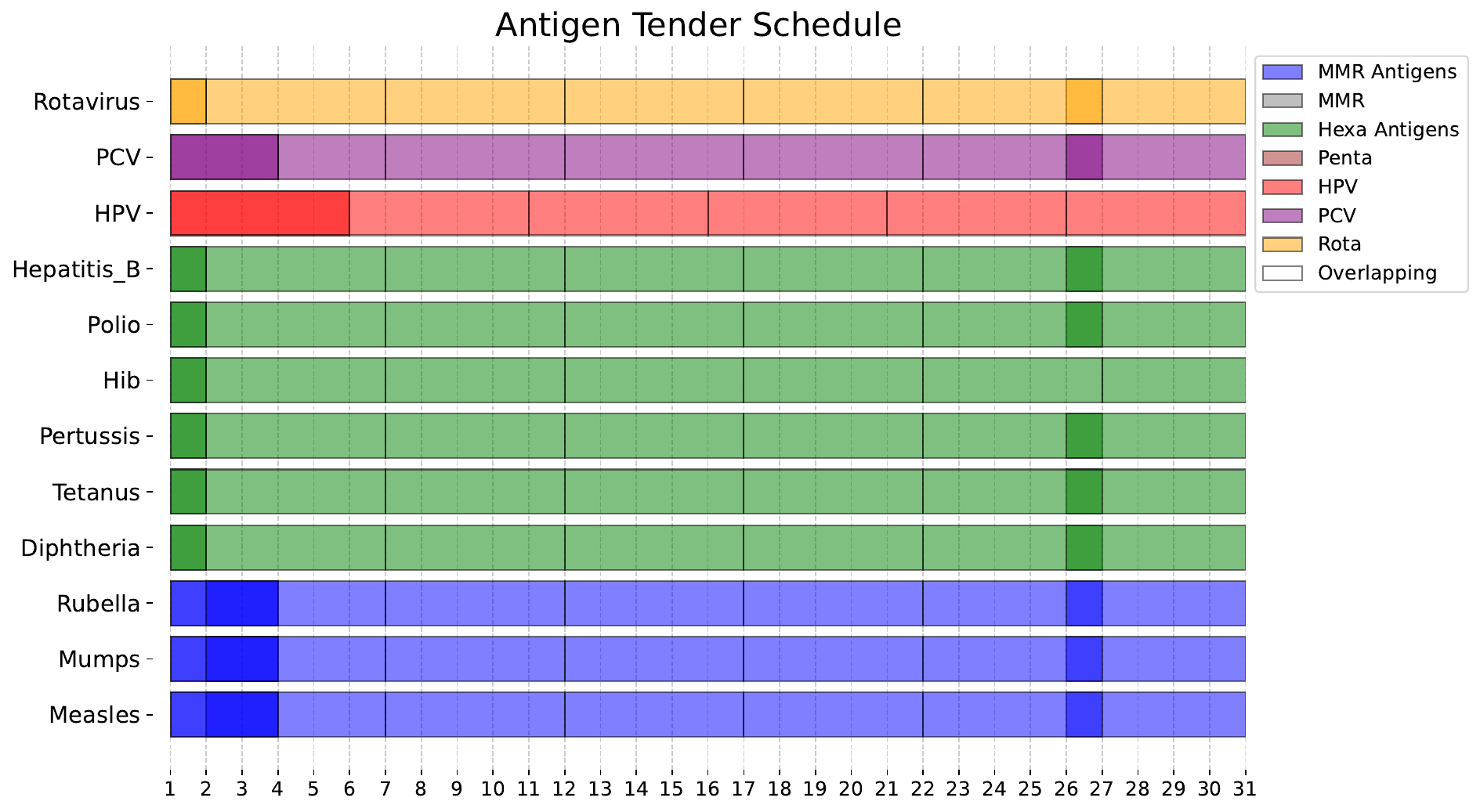}
    \caption{Tender Schedule}
    \label{fig:subfig3}
  \end{subfigure}
  \hspace{0.01\textwidth} 
  \begin{subfigure}{0.49\textwidth}
    \centering
    \includegraphics[width=\linewidth]{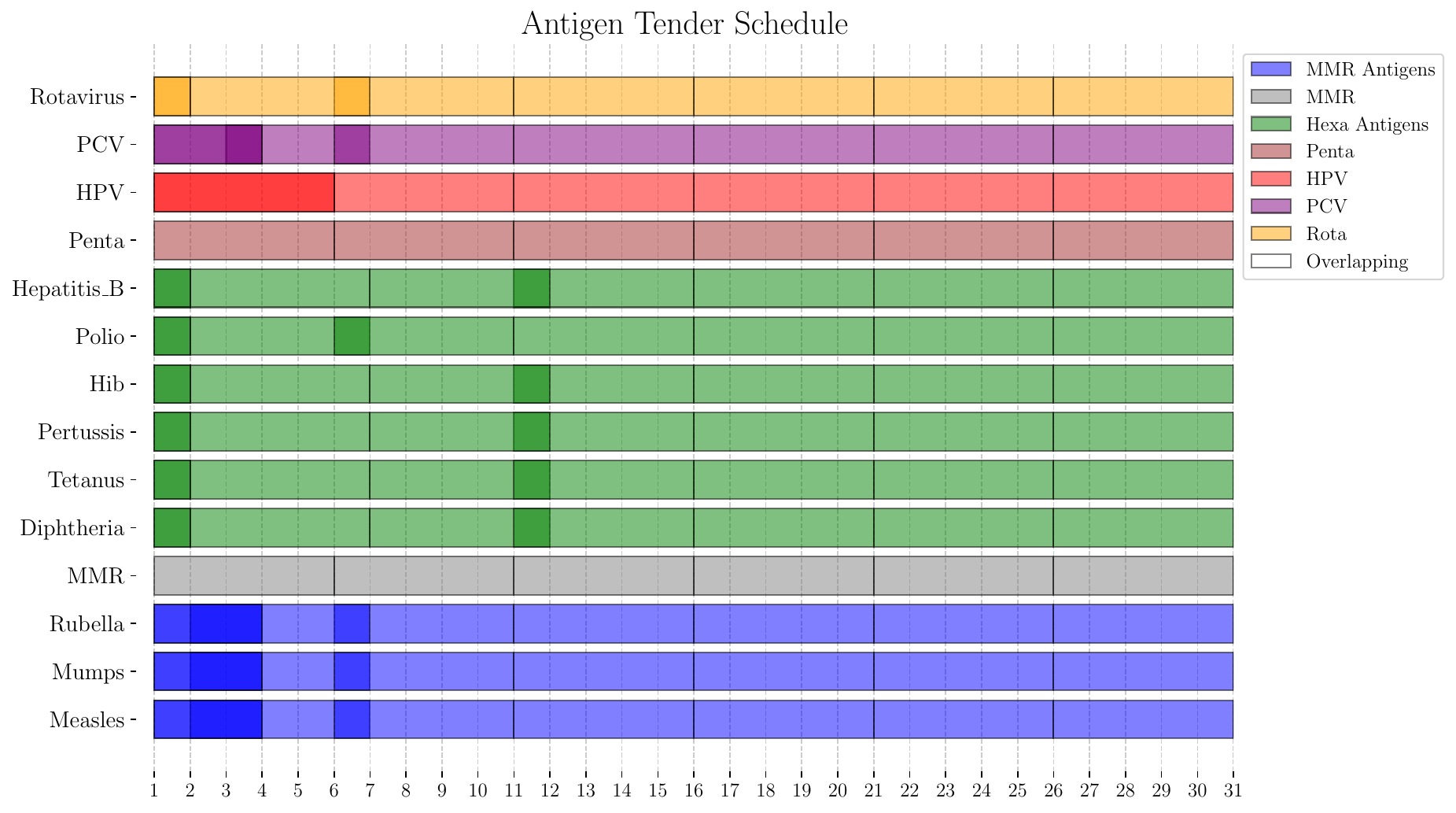}
    \caption{Tender Schedule: Vaccines as Antigens}
    \label{fig:subfig4}
  \end{subfigure}
  \caption{Tender schedules for a 30-year planning horizon}
  \label{fig:both2}
\end{figure}

Our numerical analysis determines the optimal sequencing of vaccine tenders across different planning horizons. 
Since we minimize the total discounted cost, the model tends to create long tenders. Such solutions indirectly minimize the efforts that go into organizing tenders and, thus,  the corresponding costs. As seen in Figure \ref{fig:both2}, most of the scheduled tenders have a 5-year duration. These long-term schedules show general stability, {and synchronization of common antigen sequences}. {Figure \ref{fig:subfig3} establishes a tender schedule over a 30-year period, assuming that demand is met via monovalent vaccines only. Figure \ref{fig:subfig4} shows the schedule resulting from the use of monovalent and combination vaccines. In this case, we notice that tenders still occur for all antigens. However, demand is met through the combination vaccines, where MMR and Penta account for 97\% of vaccines produced, for their corresponding antigen sets. Similarly, we notice that antigens bundled in common combination vaccines, like penta and MMR, generally stay in sequence over time. As the various color bands in Figures \ref{fig:both2} denote, antigens do not generally deviate from such bundles.}

In examining model robustness, we utilize cosine similarity, a measure to estimate the similarity between two non-zero vectors, to compare schedules. The tender schedules are transformed into vectors representing ones for a scheduled tender, for all tender durations and antigens, zeros otherwise. In our experiment, we examine varying supply and demand, and analyze 12 tender schedules. We highlight the similarity between tender schedules, with varying conditions as previously described, and demonstrate the robustness of our model. The schedules in Figures \ref{fig:subfig3} and \ref{fig:subfig4} attain a cosine similarity of 0.45, showing moderate similarity between schedules with varying antigens (single versus combination vaccine). When combination vaccines as antigens are considered in the 12 experiments, we observe cosine similarities between 0.68 and 0.79. Figure \ref{fig:subfig4} represents one such experiment.

\begin{figure}[h]
  \centering
  \begin{subfigure}{0.49\textwidth}
    \centering
    \includegraphics[width=\linewidth]{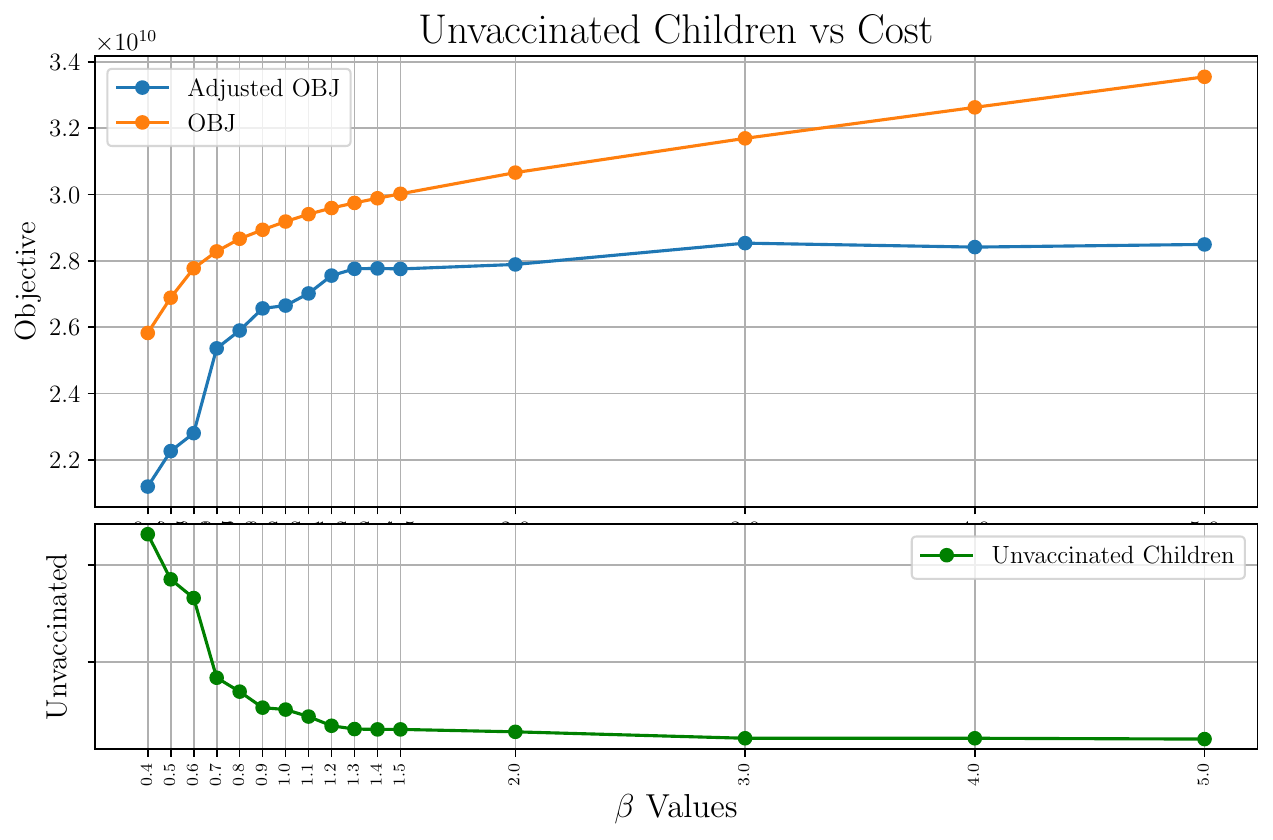}
    \caption{Unvaccinated children penalty ($\beta$)}
    \label{fig:subfig1}
  \end{subfigure}
  \hspace{0.01\textwidth} 
  \begin{subfigure}{0.49\textwidth}
    \centering
    \includegraphics[width=\linewidth]{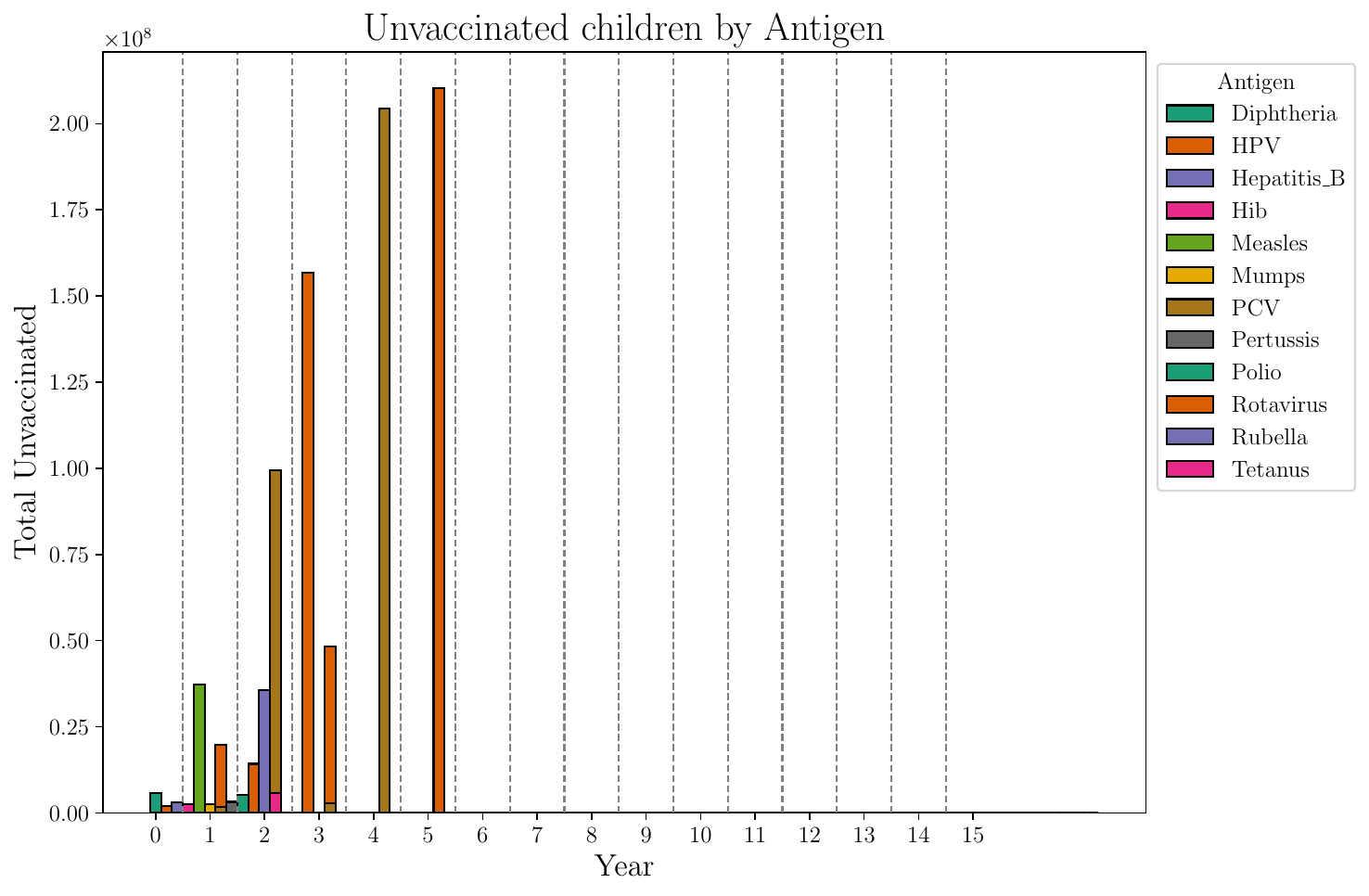}
    \caption{Unvaccinated children when $\beta$=6}
    \label{fig:subfig2}
  \end{subfigure}
  \caption{Analysis of unvaccinated children}
  \label{fig:both1}
\end{figure}

{In the objective function of our model, we include a term that penalizes the overall number of unvaccinated children.} We consider the total number of unvaccinated children as the social impact of the proposed tender schedule. The size of the penalty term $\beta$ in the objective determines the number of unvaccinated children. When this cost is (perceived to be) higher than the cost of a vaccine,  the number of unvaccinated children is lowest. Figure \ref{fig:subfig1} shows that the number of unvaccinated children decreases as the value of $\beta$ increases from \$0.4 to \$1.5. Beyond this point, increasing the value of the penalty does not account for a significant change in the number of children who are not vaccinated.

Figure \ref{fig:subfig2} presents the total number of unvaccinated children over time for $\beta = 6.$ On-going tenders and commitments, and initial inventory level drive the number of unvaccinated children during the first five years. When the value of $\beta$ is higher than the expected cost of a vaccine (as is the case with $\beta = 6$), then the number of unvaccinated children is lowest beyond the initial five years. {Unvaccinated children for Rotavirus and PCV, shown in years 4 and 5 in Figure \ref{fig:subfig2}, represent the initial conditions for this experiment. However, given model conditions, we can assess a high enough penalty for unvaccinated children to meet antigen demand by manipulating the parameter $\beta$. From this analysis, we can assess that, even with changing conditions in a hypothetically coordinated UNICEF market, we can adjust the penalty for unvaccinated children to account for these changes.}

\section{Conclusion}

In this study, we investigate the potential benefits of long-term planning for the hypothetically coordinated vaccine market (HCVM) by employing dynamic scheduling of vaccine tenders. Our approach relies on a deterministic mixed integer model to reduce the overall discounted cost of organizing tenders, manufacturing vaccines, and the cost of unvaccinated children. The proposed model facilitates the generation of vaccine tender schedules for various durations. We analyze these schedules' general stability, synchronization, and similarity and observe that the schedules generally align by antigen groups even when combination vaccines are added to the problem. The model also determines the volumes manufacturers should commit to and optimal tender lengths to satisfy the hypothetically coordinated UNICEF market demand. Future research endeavors will extend this model into a multi-stage stochastic optimization framework to account for uncertainties in demand. This extension will provide a more robust and flexible approach to vaccine tender scheduling. We anticipate that such enhancements will contribute to a more comprehensive understanding of the HCVM dynamics.

\section*{Acknowledgements}
Funding: This work was partially supported by The Bill \& Melinda Gates Foundation, Seattle, WA [grant number INV-044813].

\bibliography{bibliography}

\end{document}